\documentclass[12pt]{article}
\usepackage{graphicx}
\usepackage{amssymb}
\usepackage{amsmath}

\newtheorem{remark}{Remark}
\newtheorem{theorem}{Theorem}

\newtheorem{proposition}{Proposition}

\newcommand{\halmos}{\vspace{3mm} \hfill \mbox{$\Box$}\\[2mm]}

\begin{document}

\title{On distance in total variation between image measures}
\vspace{15pt}

\author{Youri Davydov \\
{\small Universit\'e Lille 1, Laboratoire Paul Painlev\'e }\\
{\small and}\\
{\small Saint Petersbourg State university}
}

\date{}

\vspace{20pt}

\maketitle

Abstract

We are interested in the estimation of the distance in total variation
$$
\Delta := \|P_{f(X)} - P_{g(X)}\|_{\mathrm var}
$$
between distributions of random variables $f(X)$ and $g(X)$ in terms of pro\-ximity of $f$ and $g.$
We propose a simple general method of estimating $\Delta$. For Gaussian  and trigonometrical polynomials it gives an asymptotically optimal result (when the degree tends to $\infty$).

MSC

Primary 60E05

secondary 60E15, 60A10


\section{Introduction}
 
Let $X$ be a  random vector with values in $\mathbb R^d$ having an absolutely continuous distribution $P,$ and  $f,\;g$ be two measurable functions from $\mathbb R^d$ to  $\mathbb R^1.$ We are interested in the estimation of the distance in total variation
$$
\Delta (f,g):= \|P_{f(X)} - P_{g(X)}\|_{\mathrm var}
$$
between distributions of random variables $f(X)$ and $g(X)$ in terms of pro\-ximity of $f$ and $g.$
This problem has applications in different fields of probability theory. The most interesting  example may be is the case where $d=1,\;X$ is a  Gaussian r.v. with mean $a$ and variance $\sigma^2,$ and $f,\;g$ are two polynomials of degree $m$:
$$
f(x) = \sum_0^m a_kx^{m-k},\;\;\; g(x) = \sum_0^m b_kx^{m-k},
$$
$a_k \neq 0.$


The result by Yu. Davydov and G. Martynova   (1987)  says that there exists a constant
$C$ depending only on $\;m,  a, \sigma\;$ such that
\begin{equation}\label{DM}
\Delta \leq C|a_0|^{-\frac{1}{m}}\delta^{\frac{1}{m}},
\end{equation}
where $\delta = \max_{0\leq k \leq m}|a_k - b_k|.$

The importance of this case is explained by
strong relations with the estimation of total variation distance between distributions of multiple Wiener integrals. Namely, from (\ref{DM}) it follows (for details see 
(\cite{[DM]}))
\begin{equation}\label{DM1}
\Delta({I_m(f)},\; {I_m(g)})\leq
C\|f-g\|^\frac{1}{m}_{L^2(T)},
\end{equation}
where $I_m(f), I_m(g)$ are two $m$-multiple Wiener-Ito integrals; the constant $C$ depends only on $m$ and $f$.

Below, in section 3.1, we propose some explanation of how the estimate (\ref{DM1}) could be deduced from (\ref{DM}).

In work (\cite{[P]}) an attempt to obtain an estimate for $\Delta({I_m(f)},\; {I_m(g)})$ by means of methods of stochastic analysis has been made, but it gives an order $\frac{1}{2m}$, which is significantly weaker.


When our article had been already sent for the press, we have learned about a preprint (\cite{[BKZ]}) which contains a number of the deep results connected with this problem.
In particular, it is shown that the density of distribution of any non-constant Gaussian polynomial of degree $m$ always belongs to the Nikol'ski-Besov class $B^{\frac{1}{m}}(\mathbb R^1)$,
and in the one-dimensional case the estimate (\ref{DM}) is proved with logarithmic factor.

The aim of the present work is to propose a simple general method of estimating $\Delta(f,g)$.
For completely different reasons we independently
arrived to the use of condition type (\ref{f1})
and showed (see Th. 1) that having this condition
(in arbitrary dimension) with the exponent $\alpha$,
we obtain for $\Delta(f,g)$ the order $\frac{\alpha}{\alpha +1}.$ 
In combination with the aforementioned result from (\cite{[BKZ]}) it follows from our Th.1 that for Gaussian polynomials in any dimension
$$
\Delta(f,g) = O\left(\|f-g\|^\frac{1}{m+1}_{L^2(T)}\right),
$$
 which will still be asymptotically optimal (when the degree $m$ tends to $\infty$).

 As a second example we consider the case where $f$ and $g$ are trigonometrical polynomials. Here also our method gives an asymptotically optimal estimate. 

\section{Results}
     We use the notation $\cal P$ for the distribution of $X$ and  $\|\cdot\|_1$ for the norm in the space
${{\mathbf L}^1(d{\cal P})}$ of integrable functions with respect to the measure $\cal P.$

Recall that for a signed measure $\mu$ its total variation is defined by
$$
\|\mu\|_{\mathrm var} = \sup_{(A_i)}\sum_i|\mu(A_i)|,
$$
where the supremum is taken over all finite measurable partitions $(A_i)$ of the space. If $\mu$ has a density $m$ with respect to some non negative measure $\nu$ then
$$
\|\mu\|_{\mathrm var} = \int|m|d\nu.
$$
\begin{theorem} \hspace{-10pt}
{\huge .}\;Suppose that for some $\alpha > 0$
\begin{equation}\label{f1}
\|P_{f(X)} - P_{f(X)+u}\|_{\mathrm var} \leq C_f|u|^\alpha,
\end{equation}
and
\begin{equation}\label{g1}
\|P_{g(X)} - P_{g(X)+u}\|_{\mathrm var} \leq C_g|u|^\alpha.
\end{equation}
Then
\begin{equation}\label{fg}
\|P_{f(X)} - P_{g(X)}\|_{\mathrm var} \leq C\|f-g\|_1^{\frac{\alpha}{\alpha+1}},
\end{equation}
where $C= (C_f+C_g)^{\frac{1}{\alpha +1}}
(E|\nu|^\alpha  \;+\;\sqrt{\frac{\pi}{2}}),$ and $\nu$ is a standard Gaussian r.v.
\end{theorem}

\begin{remark}\hspace{-10pt}
{\huge .}\;
It is known that $E|\nu|^\alpha = \frac{2^{\frac{\alpha}{2}}{\mathbf \Gamma}(\frac{\alpha+1}{2})}{\sqrt{\pi}},$ where
$\mathbf\Gamma$ is the Gamma function.
\end{remark}

\begin{remark}\hspace{-10pt}
{\huge .}\;
As we always have $\|P_{f(X)} - P_{g(X)}\|_{\mathrm var} \leq 2,$ one can replace the expression in the right part of (\ref{fg})
by $\min\{2,\;\; C\|f-g\|_1^{\frac{\alpha}{\alpha+1}}\}.$
\end{remark}

\begin{remark}\hspace{-10pt}
{\huge .}\;
In the case when $P_{f(X)}$ has a density $\;p$ with respect to the Lebesgue measure the condition (\ref{f1})
means that  $p$ belongs to the so-called Nikol'ski-Besov space $B_1^{\alpha}(\mathbb R^1)$ (for details see (\cite{[BKZ]}).
\end{remark}

{\bf Proof.}
Let $\nu$ be a standard Gaussian r.v. independent of $X$ and $\xi = \sigma\nu$ where $\sigma$ is a positive number,
its exact value  will be chosen later. 

We have
$$
\|P_{f(X)} - P_{g(X)}\|_{\mathrm var} \leq \delta_1 + \delta_2 +\delta _3,
$$
where
$$
\delta_1 = \|P_{f(X)} - P_{f(X)+\xi}\|_{\mathrm var},
$$
$$
\delta_2 = \|P_{g(X)} - P_{g(X)+\xi}\|_{\mathrm var},
$$
$$
\delta_3 = \|P_{f(X)+\xi} - P_{g(X)+\xi}\|_{\mathrm var}.
$$
We find using (\ref{f1})
$$
\delta_1 = \|P_{f(X)} - P_{f(X)+\xi}\|_{\mathrm var} = 
\int_{\mathbb R}\|P_{f(X)} - P_{f(X)+u}\|_{\mathrm var}P_\xi(du) \leq
$$
\begin{equation} \label{delta1}
\leq C_f\int_{\mathbb R}|u|^\alpha P_\xi(du) = C_f E|\xi|^\alpha = C_f\sigma^\alpha E|\nu|^\alpha.
\end{equation}
Similarly
\begin{equation} \label{delta2}
\delta_2  \leq C_g\sigma^\alpha E|\nu|^\alpha.
\end{equation}
Consider now $\delta_3.$ Denoting $\tilde P,\;\tilde Q$ distributions in ${\mathbb R}^2$ of random vectors
$(X,\,f(X)+\xi)$ and $(X,\,g(X)+\xi),$  we remark that
$$
P_{f(X)+\xi} =\tilde Ph^{-1},\;\;P_{g(X)+\xi} =\tilde Qh^{-1},
$$
where $h : {\mathbb R}^2 \rightarrow  {\mathbb R}, \;\; h(x,y) = y.$
Therefore,
$$
\delta_3 =  \|P_{f(X)+\xi} - P_{g(X)+\xi}\|_{\mathrm var} \leq \|\tilde P - \tilde Q\|_{\mathrm var}.
$$
It is easy to see that
\begin{equation} \label{delta3}
\|\tilde P - \tilde Q\|_{\mathrm var}  \leq \int_{\mathbb R}\|P_{f(x)+\xi} - P_{g(x)+\xi}\|_{\mathrm var}P_X(dx).
\end{equation}
As the distributions $P_{f(x)+\xi}$ and $P_{g(x)+\xi}$ are Gaussian with the same variance $\sigma^2$ and with 
mean values
differing by $|f(x)-g(x)|,$  we have
$$
 \|P_{f(x)+\xi} - P_{g(x)+\xi}\|_{\mathrm var} \leq \frac{2}{\sigma\sqrt{2\pi}}|f(x)-g(x)| .
$$
Hence, it follows  from (\ref{delta3}) that
\begin{equation} \label{delta3a}
\|\tilde P - \tilde Q\|_{\mathrm var}  \leq \frac{2}{\sigma\sqrt{2\pi}}\|f-g\|_1.
\end{equation}
Gathering estimates (\ref{delta1}),  (\ref{delta2}) and  (\ref{delta3a}), we get
$$
\|P_{f(X)} - P_{g(X)}\|_{\mathrm var} \leq  (C_f+C_g)\sigma^\alpha E|\nu|^\alpha +
\frac{2}{\sigma\sqrt{2\pi}}\|f-g\|_1.
$$
Taking $\sigma = \{(C_f + C_g)\|f-g\|_1\}^\frac{1}{1+\alpha},$
we find the final result.
\halmos

 Suppose now that the dimension $d=1$ and consider some sufficient conditions for the relations of type (\ref{f1}).
Remarking that using the notation $f_u(t)$ for $f(t-u)$, and $\cal P$ for the distribution of $X$, we can rewrite
the value $\delta(u)= \|P_{f(X)} - P_{f(X)+u}\|_{\mathrm var}$ in the equivalent form:
$$
\delta(u)= \|{\cal P}{f^{-1}} - {\cal P}{f_u^{-1}}\|_{\mathrm var}.
$$

Below we will also use this notation in the case where $\cal P$ is finite but not necessarily a probability measure. 
\begin{proposition}\label{p1}\hspace{-10pt}
{\huge .}\;
Let $f$ be a convex strictly increasing  function defined on the interval $[a,b]$ and such that for some $m >0,\; K>0,$
\begin{equation}\label{asymp1}
f(x)-f(a) \sim  K(x-a)^m,\;\; x\downarrow a.
\end{equation}
Let ${\cal P} = \lambda\ ,\;\; \lambda$ being Lebesgue measure. 

Then 
\begin{equation}\label{asymp2}
\delta(u) \leq 2C_fu^\frac{1}{m},\;\;\; u\geq 0,
\end{equation}
where 
\begin{equation}\label{Cf}
C_f = K^{-\frac{1}{m}}\sup_{f(a)<x<f(b)}\left\{\left|\frac{f^{-1}(x)-a}{(x- f(a))^{\frac{1}{m}}}\right|\right\}.
\end{equation}
\end{proposition}

\begin{remark}\hspace{-10pt}
{\huge .}\;
It is clear that similar estimates (with evident changes) are available if we replace "convex" by "concave" and (or) "increasing" by "decreasing".
\end{remark}

{\bf Proof.} 
First of all remark that by
(\ref{asymp1}),  $f^{-1}(f(a)+u)-a \sim K^{-\frac{1}{m}}u^{\frac{1}{m}},$ when $u\rightarrow 0$, which shows that the constant $C_f$ is finite.
As $f^\prime(t) > 0$ for all $t$, the measure $\lambda f^{-1}$ has a density
$$
h(t) = \frac{1}{f'(f^{-1}(t))}{\mathbf 1}_{[f(a), f(b)]}(t)
$$
which is decreasing.

Therefore for $u \in [f(a), f(b)]$
$$
\delta (u) = 2\int_{f(a)}^{f(a)+u }h(t)dt = 2(f^{-1}(f(a)+u)-a).
$$
(The first equality will be evident if we consider the epigraphs of the functions $f(t)$ and $f(t-u).$)
Again by (\ref{asymp1}),  $f^{-1}(f(a)+u)-a \sim K^{-\frac{1}{m}}u^{\frac{1}{m}},$
which gives (\ref{asymp2}). \halmos

A more general and more useful result is given by the following proposition.

\begin{proposition}\label{p2}\hspace{-10pt}
{\huge .}\;
Let $f$ be a  convex strictly increasing  function defined on $[a,b]$ and such that for some $m >0,\; K>0,$
\begin{equation}\label{asymp3}
f(x)-f(a) \sim K (x-a)^m,\;\; x\downarrow a.
\end{equation}
Let $\cal P $ be a finite  measure on $[a,b]$ having a density $p$ which satisfies the Lip\-schitz condition:
$$
|p(x) - p(y)| \leq L|x-y|,\;\;\;\; \forall \;\;x,y \in [a,b].
$$
Let $A=\sup_{x\in [a,b]} p(x).$
Then 
\begin{equation}\label{asymp4}
\delta(u)= \|{\cal P}f^{-1} - {\cal P}f_u^{-1}\|_{\mathrm var} \leq  [3A+L(b-a)]C_fu^\frac{1}{m},\; \;\;u\geq 0,
\end{equation}
where $C_f $ is given by (\ref{Cf}).
\end{proposition}

{\bf Proof.}
The measure ${\cal P}f^{-1}$ is absolutely continuous and its density is equal to
\begin{equation}\label{q}
q(t) = h(t)p(f^{-1}(t)),
\end{equation}
where $h(t)= \frac{1}{f'(f^{-1}(t))}{\mathbf 1}_{[f(a), f(b)]}(t)$ is the density of $\lambda f^{-1}.$
Hence
$$
\delta(u) = \int_{f(a)}^{f(b)+u}|q(t)-q(t-u)|dt = I_1+I_2+I_3,
$$
where
$$
I_1 =  \int_{f(a)}^{f(a)+u}|q(t)-q(t-u)|dt,
$$
$$
 I_2= \int_{f(a)+u}^{f(b)}|q(t)-q(t-u)|dt,
$$
$$
 I_3= \int_{f(b)}^{f(b)+u}|q(t)-q(t-u)|dt.
$$

Consider $I_1.$ Since $p$ is bounded and $q(t-u) = 0$ for $t\leq u$, we have as before
$$
I_1 \leq A\int_{f(a)}^{{f(a)}+u}h(t)dt = A\lambda([a, f^{-1}(f(a)+u)]) \leq AC_fu^\frac{1}{m}.
$$

Since $h$ is decreasing, we get similarly
$$
I_3 \leq A\int_{f(b)}^{f(b)+u}h(t)dt \leq  A\int_{f(a)}^{{f(a)}+u}h(t)dt \leq AC_fu^\frac{1}{m}.
$$

By the triangle inequality
$$
I_2 \leq J_1 +J_2,
$$
where
$$
J_1 = \int_{f(a)+u}^{f(b)}h(t)|p(f^{-1}(t))-p(f^{-1}(t-u))|dt,
$$
$$
J_2 = \int_{f(a)+u}^{f(b)}p(f^{-1}(t))|h(t)-h(t-u)|.
$$
Since $p$ is Lipschitz, 
$$
J_1 \leq L \int_{f(a)+u}^{f(b)}h(t)|f^{-1}(t)-f^{-1}(t-u)|dt.
$$
As $f$ is convex and increasing, $f^{-1}$ is concave and increasing. Therefore
$$
|f^{-1}(t)-f^{-1}(t-u)| \leq |f^{-1}(f(a)+u)-a|.
$$
Hence, using that $\int_{f(a)+u}^{f(b)}h(t)dt \leq b-a,$ we get
$$
J_1 \leq L(b-a)C_fu^\frac{1}{m}.
$$
It is clear that 
$$
J_2 \leq A\int_{f(a)}^{f(b)}|h(t)-h(t-u)|,
$$
which is less than or equal to $AC_fu^\frac{1}{m}$ by Proposition 1.

Finally, gathering all previous estimations, we have
$$
\delta(u) \leq [3A+L(b-a)]C_fu^\frac{1}{m}.
$$
\halmos

\section {Gaussian polynomials}
As a first example of application we consider the case where $f,\;g$ are two polynomials of degree $m$
of $d$ variables and $\cal P$ is a standard Gaussian measure in $\mathbb R^d.$

Let $\|\nabla f\|_\star ^2 = $
$\sup_e\int_{\mathbb R^d}|\partial_e f|^2 d{\cal P},$
where $\partial_e f$ is the derivative of $f$ in the direction $e \in S_{d-1}.$
\begin{theorem}\label{DB}\hspace{-10pt}
{\huge .}\;
If $f,g$ are non-constant, then 
there exists a constant $C>0$ depending only on 
$m,\; \|\nabla f\|_\star, \;\|\nabla g\|_\star,$ such that 
\begin{equation}\label{DB1}
\Delta(f,g) \leq C\|f-g\|_1^{\frac{1}{m+1}}.
\end{equation} 
\end{theorem}

{\bf Proof.} From Th. 5.7 of \cite{[BKZ]} it follows that the conditions (\ref{f1}), (\ref{g1}) are fulfilled with $\alpha = \frac{1}{m}.$ Therefore by Th. 1 we get (\ref{DB1}). \halmos

Let us consider the one-dimensional case. Then
$$
f(x) = \sum_0^m a_kx^{m-k},\;\;\; g(x) = \sum_0^m b_kx^{m-k},
$$
$a_k \neq 0,$ 
and from  (\ref{DB1}) we deduce  the estimation
$$
\Delta(f,g) \leq C\delta^{\frac{1}{m+1}},
$$
 where $\delta = \max_{0\leq k \leq m}|a_k - b_k|.$

The order $\frac{1}{m+1}$ is worse than one in (\ref{DM}) but asymptotically (when $m \rightarrow \infty$) they are equal.

Due to the importance of condition type (\ref{f1})  it seems  reasonable to present here its elementary proof.

Let $x_1, x_2,\ldots, x_n,$ be the ordered set of 
all the roots of the derivatives $f'$ and $f^{(2)}.$ It is clear that $n\leq 2m-3.$
On each segment $\Delta_k = [x_k,x_{k+1}]$
the function $f$ is convex or concave and $f'$ can be equal to zero not more than in one of the ends 
of the segment. It means that $f$ on $\Delta_k$ satisfies condition (\ref{Cf}) for
some $m_k\leq m.$ Denote ${\cal P}_k = {\cal P}_{\Delta_k}$ the restriction of ${\cal P}$ on $\Delta_k.$ Then,
by Proposition 2,
\begin{equation}
\delta_k(u) = \|{\cal P}_kf^{-1}-{\cal P}_kf_u^{-1}\|_{\mathrm var}\leq [3A_k + L_k(x_{k+1}-x_k)]C_{f,k}u^ \frac{1}{m},
\end{equation}
where $A_k= \sup_{x\in \Delta_k}p(x),\;\;\; L_k =  \sup_{x\in \Delta_k}p'(x),$ and  $C_{f,k} $ is defined by (\ref{Cf})
with $a=x_k, b=x_{k+1}$ and $d$ depending on $ \Delta_k.$

Summing these estimates, we find
\begin{equation}\label{rho-0}
\rho_0(u): = \|{\cal P}_{[x_0,x_n]}f^{-1}-{\cal P}_{[x_0,x_n]}f_u^{-1}\|_{\mathrm var}\leq C_1u^ \frac{1}{m},
\end{equation}
where 
$$
C_1 =\sum_{k=0}^{n-1}[3A_k +L_k(x_{k+1}-x_k)]C_{f,k}.
$$
To estimate
$$
\rho_+(u): = \|{\cal P}_{[x_n,\infty)}f^{-1}-{\cal P}_{[x_n,\infty)}f_u^{-1}\|_{\mathrm var}
$$
we represent  $[x_n,\infty)$ as the union of segments :\\
$[x_n,\infty)=\cup_{j=0}^{\infty}\Delta_j,
[x_n+j, x_n+j+1].$
Similarly to before, we get
\begin{equation}\label{rho-+}
\rho_+(u) \leq C_2u^ \frac{1}{m},
\end{equation}
where now
$$
C_2 =\sum_{k=0}^{\infty}[3A_k +L_k]C_{f,k}.
$$
Since $f$ is convex on $[x_n,\infty),\;\;\; C_{f,k} \leq C_{f,0}.$ 
The series $\sum_{k=0}^{\infty}[3A_k +L_k]$ is convergent because $p$ is Gaussian density. Therefore the constant
$C_2$ is finite.

Applying similar arguments to the estimation of 
$$
\rho_-(u): = \|{\cal P}_{(-\infty,x_0]}f^{-1}-{\cal P}_{(-\infty,x_0]}f_u^{-1}\|_{\mathrm var},
$$
we see that 
$$
\rho_-(u) \leq C_3u^ \frac{1}{m}
$$
for some $C_3 <\infty.$
This inequality together with (\ref{rho-0}) and (\ref{rho-+}) gives the final result:
 the condition (\ref{f1}) is fulfilled for $f$ with $\alpha = \frac{1}{m}.$
 
 \subsection{Multiple integrals}
 
 Let $W$ be random Gaussian orthogonal measure corresponding to the Le\-besgue measure $\lambda$ on 
 $\mathbb R^1,\;\; EW(A) = 0,  EW(A)W(B) = \lambda(A\cap B).$ Let $H_n$ be the space of functions $f: \mathbb R^n \rightarrow 
 \mathbb R^1$ which are square integrable with respect to $\lambda^n$ and are invariant under all permutations of coordinates. For such a function the multiple integral
 $$
 I_n(f) = \int_{\mathbb R^n}\,f(x_1,\ldots,x_n)W(dx_1)\ldots
 W(dx_n)
 $$
 is well defined (see for details \cite{[M1]}, \cite{[D]}).
 
 Let $P$ be the distribution of $W$ in the space 
 ${\mathbb S}=(\mathbb R^{\cal A},{\cal B}^{\cal A}),$
 where ${\cal A}= \{A\in {\cal B}^1\;|\; \lambda (A) < \infty\}.$
 
 The measure $P$ is  Gaussian  and its admissible shifts $\nu = \nu_h$ are exactly  the measures which are absolutely continuous  with respect to $\lambda$ (see Prop. 2, \cite{[D]}) and such that
$\nu_h(A) = \int_A hd\lambda,\;\;\;h\in L^2(d\lambda).$
If $\Gamma$ is a partition of ${\mathbb S}$ composed by the lines $\{l_\varkappa=\varkappa+c\nu_h,\; c\in \mathbb R^1\}$
parallel to $\nu_h,$ then the conditional distributions $(P_\varkappa)$ for $P$ on these lines
will be Gaussian with the mean value 
$a_h=-\|h\|^2_{H_1} \int h d\varkappa$ and the variance
$\sigma^2_h=\|h\|^{-2}_{H_1} $ (see Prop. 3, \cite{[D]}).

The integral $I_n(f)$ can be considered as a measurable functional
$$
I_n(f)= F(\varkappa) = \int f d\varkappa
$$
and its restriction onto $l_\varkappa$ is a polynomial
of the degree $n$:
$$
F_\varkappa(c) = F(\varkappa+c\nu_h) = 
c^n\int fd_n\nu_h + \sum_{m=0}^{n-1}\xi_m c^m,
$$
where $\xi_m$  are some functions on $\varkappa$ and $\nu_h.$
In \cite{[D]} it is shown that we can choose $\nu_h$ in such a way that $\int fd_n\nu_h\neq 0.$
Hence $F_\varkappa$ is a polynomial of the degree $n$
and the measure $P_{I_n(f)}$ can be represented as a mixture of distributions of one-dimensional Gaussian polynomials
$$
P_{I_n(f)} = \int_{{\mathbb S}/\Gamma}P_\varkappa
F_\varkappa^{-1}\;P_\Gamma(d\varkappa),
$$
where $P_\Gamma$ is the factor-measure.

Similarly,
$$
P_{I_n(g)} = \int_{{\mathbb S}/\Gamma}P_\varkappa
G_\varkappa^{-1}\;P_\Gamma(d\varkappa),
$$
where $G_\varkappa$ is the restriction of $I_n(g)$ onto 
$l_\varkappa.$ 

Therefore
\begin{equation}\label{MI}
\|P_{I_n(f)}- P_{I_n(g)}\|\leq 
\int_{{\mathbb S}/\Gamma}\|P_\varkappa
F_\varkappa^{-1}- P_\varkappa
G_\varkappa^{-1}\|\;P_\Gamma(d\varkappa).
\end{equation}

Without loss of generality we can suppose additionally that $h$ is conti\-nu\-ous. Then we can identify the factor-space ${\mathbb S}/\Gamma$ with the subspace\\
$\{\varkappa\in {\mathbb S}\,|\,\int h\,d\varkappa =0\}.$ At the same time conditional measures $P_\varkappa$ will be Gaussian with parameters $(0, \sigma^2_h)$ which don't depend on $\varkappa.$
Hence from (\ref{MI}) and one-dimensional estimate (\ref{DM}) we directly deduce (\ref{DM1}).

\section { Trigonometrical polynomials}
As a second example we consider the case where $f$ and $g$ are two trigonometrical polynomials:
$$
f= \sum_{k=0}^n(a_k\cos kx + b_k\sin kx),\;\;\; g=  \sum_{k=0}^n(c_k\cos kx + d_k\sin kx).
$$
Like before, we suppose that ${\cal P}$ is a standard Gaussian distribution.

It is clear that the exponent $\alpha$ in (\ref{f1}) depends on the number $\varkappa$ of zero derivatives at fixed points
of the function $f.$  Let us show that in general that number cannot be more than $2n-1.$

Consider the polynomial $f.$ Without loss of generality we can and do suppose that $a_0=0$ and $x=0.$
The assertion $f^{(l)}(0)=0$ for $l = 1,\ldots,2m$ is equivalent to the statement that the system of $2m$ linear equations  
(with respect to unknowns $a_k$ and $b_k,\,k =1,\ldots,n$)
$$
\begin{array}{cc}  
\begin{cases} 
\sum_1^n kb_k=0 \\ \sum_1^n k^3b_k=0 \\ \cdots\\ \sum_1^n k^{2m-1}b_k=0 
\end{cases} 
&
\;\;\;\;\;\begin{cases} 
\sum_1^n k^2a_k=0 \\ \sum_1^n k^4a_k=0 \\ \cdots\\ \sum_1^n k^{2m}a_k=0 
\end{cases} 
\end{array}
$$
has a non-trivial solution.

For $m=n$ the determinant $\Delta $ of this system satisfies the following relation
$$
\Delta = (n!)^3W^2(1,2^2,3^2,\ldots,n^2),
$$
where $W(x_1, \ldots,x_n)$ is the Vandermonde determinant.

Hence $\Delta \neq 0$ and therefore our system can have non-trivial solution only if $2m\leq 2n-1.$ It means that in general  $\varkappa \leq 2n-1.$ 
The case $m$ is odd gives the same born.
Now, arguments similar to ones used in the previous section show that the conditions 
(\ref{f1}), (\ref{g1}) hold with $\alpha = \frac{1}{2n}.$ By Th. 1 we get in this case the following estimation
$$
 \|{\cal P}{f^{-1}} - {\cal P}{g^{-1}}\|_{\mathrm var} \leq C\|f-g\|_\infty^{\frac{1}{2n+1}}.
$$

\section { Concluding remarks}
1. In multi-dimensional setting in the class of all polynomials the order $\frac{1}{m+1}$
in the estimate (\ref{DB1}) is asymptotically the best possible. At the same time the example of the polynomial $f(x_1,\ldots,x_d)= (x_1^2 +\cdots+x_d^2)^m$
shows that (\ref{DB1}) is fulfilled with the exponent
$\min\{1, \frac{d}{m}\}.$ It would be interesting to describe precisely the sub classes of polynomials which provide intermediate orders.

\vspace{7pt}

\noindent
2. The proof of (\ref{DM}) in (\cite{[DM]})
is strongly based on  the particular properties  of  usual polynomials and cannot be applied even in the case of trigonometrical polynomials.  It would be interesting to find a general approach which allows to reach optimal estimates.
\vspace{7pt}

\noindent
3. It would be also interesting to find sufficient conditions for the application of our Th. 1 to analytic functions $f$ and $g.$
\vspace{10pt}

{\bf  Acknowledgments}
\vspace{7pt}

I am very grateful to the anonymous referee for the reference to important work [1] and
for the competent remarks which have allowed to improve significantly our article.

\section { References}

\bibliographystyle{plain}

\begin{enumerate}


\bibitem {[BKZ]}   Bogachev V. I.,  Kosov E. D. and  Zelenov G. I., 2016. Fracti\-onal smoothness of distributions of polynomials and a fractional analog of Hardy-Landau-Littlewood inequality. 
 preprint 
arXiv:1602.05207.

\bibitem  {[DLS]} Davydov Yu., Lifshits M. A. and Smorodina N. V., 1998. Local properties of distributions of stochastic functionals. New-York : AMS eds.,
184p.

\bibitem  {[D]} Davydov Yu., 1990. On distributions of Wiener-Ito multiple stochastic
integrals. Theor. Probab. Appl., v. 35, 1, 27-37.

\bibitem  {[DM]}  Davydov Y. A.,  Martynova G. V., 1987. Limit behaviour of multiple stochastic integral. Statistics and control of random process. Preila, Nauka,
Moscow, 55--57 (in Russian).

\bibitem  {[M]}  Martynova G. V., 1987. PhD thesis, Saint Petersbourg state university.

\bibitem  {[M1]} Major P.,  Multiple Wiener-Ito Integrals.— Lect. Notes Math., 1981. B. 849, S. 127.

\bibitem  {[P]}  Nourdin I., Poly G., 2013. Convergence in total variation on Wiener chaos. Stoch. Proc. Appl. 123, 651--674.

\end{enumerate}

\end{document}